\documentclass{amsart}
\usepackage{amsmath,amsthm,amssymb}

\makeatletter
    
    \@addtoreset{equation}{section}
  \makeatother

\newtheorem{definition}{Definition}[section]

\newtheorem{theorem}[definition]{Theorem}

\pagebreak

\title[Semilinear damped wave equation on the half-space]{
A remark on the critical exponent
for the semilinear damped wave equation
on the half-space
}
\author[Y. Wakasugi]{Yuta Wakasugi}
\address[Y. Wakasugi]{
Department of Engineering for Production and Environment,
Graduate School of Science and Engineering,
Ehime University,
3 Bunkyo-cho, Matsuyama, Ehime, 790-8577, Japan
}
\email{wakasugi.yuta.vi@ehime-u.ac.jp}
\begin{document}
\begin{abstract}
In this short notice, we prove the non-existence of global solutions to
the semilinear damped wave equation on the half-space,
and we determine the critical exponent for any space dimension.
\end{abstract}
\keywords{semilinear damped wave equation; half-space; blow-up; critical exponent}

\maketitle
\section{Introduction}

Let $n \ge 1$ be an integer and
let $\mathbb{R}^n_+$ be the $n$-dimensional half-space, namely,
\[
	\mathbb{R}^n_+ = \{ x = (x_1, \ldots, x_n) \in \mathbb{R}^n \,;\, x_n > 0 \}
	\ (n \ge 2),\quad
	\mathbb{R}_+ = (0,\infty)
	\ (n=1).
\] 
We consider the initial-boundary value problem for the
semilinear damped wave equation on the half-space:
\begin{align}%
\label{dw}
	\left\{\begin{array}{ll}
	u_{tt}-\Delta u+ u_t = |u|^p&t>0, x\in \mathbb{R}^n_+,\\
	u(t,x) = 0,&t>0, x \in \partial \mathbb{R}^n_+,\\
	u(0,x) = u_0(x), \ u_t(0,x) = u_1(x),&x\in \mathbb{R}^n_+.
	\end{array}\right.
\end{align}%
Here, $u$ is a real-valued unknown function and
$u_0, u_1$ are given initial data.

Our aim is to show the non-existence of global solutions
and determine the critical exponent for any space dimension.
Here, the {\em critical exponent} stands for the threshold of the exponent of the nonlinearity
for the global existence and the finite time blow-up of solution with small data.

For the semilinear heat equation $v_t - \Delta v = v^p$ on the whole space,
Fujita \cite{Fu66} discovered that
if $p > p_F(n) := 1+2/n$, then the unique global solution exists for every small positive initial data,
while the local solution blows up in finite time for any positive data if
$1<p<p_F(n)$.
Namely, the critical exponent of the semilinear heat equation on the whole space
is given by $p_F(n)$, which is so-called Fujita's critical exponent. 
Later on, Hayakawa \cite{Hay73} and Kobayashi, Shirao and Tanaka \cite{KoShTa77}
proved that the case $p=p_F(n)$ belongs to the blow-up region.
Moreover, the initial-boundary value problem of the semilinear heat equation
on the halved space $\mathbb{R}^k_+ \times \mathbb{R}^{n-k}_+$
was studied by \cite{LeMe89, LeMe90, Me88, Me90}
and they determined the critical exponent as
$p=p_F(n+k)$.

The critical exponent for the semilinear damped wave equation on the whole space
was studied by many authors and it is determined as $p = p_F(n)$.
We refer the reader to
\cite{ToYo01, Zh01}
and the references therein.

Ikehata \cite{Ik03DIE, Ik03JMAA, Ik04MMAS} studied the semilinear damped wave equation
on the half-space \eqref{dw} and proved that
if
$p_F(n+1) < p < \infty\ (n=1,2),\ p_F(n+1) < p \le \frac{n}{n-2} \ (n\ge 3)$,
$(u_0, u_1) \in H^1_0(\mathbb{R}^n_+) \times L^2(\mathbb{R}^n_+)$
have compact support in $\overline{\mathbb{R}^n_+}$ and
$\| \nabla u_0 \|_{L^2} + \|u_1 \|_{L^2}$ is sufficiently small,
then the problem \eqref{dw} admits a unique global solution.
When $n=1$,
Nishihara and Zhao \cite{NiZh05} proved the blow-up of solutions
when $1<p \le p_F(2)$, namely, the critical exponent of \eqref{dw}
on the half-line is determined as $p=p_F(2)$.
However, there is no blow-up result for \eqref{dw} when $n\ge 2$.

In this paper, we prove the non-existence of global classical solutions
for \eqref{dw} for all $n\ge 1$,
and we determine the critical exponent of \eqref{dw} as $p_F(n+1)$.

\begin{theorem}
Let
$1< p \le p_F(n+1) = 1+ \frac{2}{n+1}$.
We assume that the initial data satisfy
$x_n u_0, x_n u_1 \in L^1(\mathbb{R}^n_+)$
and
\begin{align}%
\label{ini}
	\int_{\mathbb{R}^n_+} x_n (u_0(x) + u_1(x))\,dx > 0
\end{align}%
(when $n=1$, we interpret $x_n = x$).
Then, there is no global classical solution to \eqref{dw}.
\end{theorem}

Our proof is based on the test function method by Zhang \cite{Zh01}.
To apply it to the half-space, we employ the technique by
Geng, Yang and Lai \cite{GeYaLa16}.
Namely, we use the test function having the form
$x_n \psi_R(t,x)$, where $\psi_R (t,x)$ is a test function supported
on the rectangle $\{ (t,x) \in [0,\infty) \times \mathbb{R}^n
\,;\, t \le R^2, |x_j| \le R \ (j=1, \ldots, n)\}$.

\section{Proof of Theorem 1.1}
We suppose that the global classical solution $u$ of the problem \eqref{dw}
exists and derive the contradiction.
Let $\psi \in C_0^{\infty}([0,\infty) \times \mathbb{R}^n_+)$ be a test function.
Using the integration by parts, we compute
\begin{align}%
\label{int}
	\int_0^{\infty} \int_{\mathbb{R}^n_+} |u|^p \psi\,dxdt
	&= \int_0^{\infty} \int_{\mathbb{R}^n_+} ( u_{tt} - \Delta u + u_t) \psi \,dxdt\\
\nonumber
	&= \int_0^{\infty} \int_{\mathbb{R}^{n-1}} \partial_{x_n} u(t,x',0) \psi (t,x',0) \,dx' dt\\
\nonumber
	&\quad+
		\int_0^{\infty}\int_{\mathbb{R}^n_+}
			u (\psi_{tt} - \Delta \psi - \psi_t ) \,dxdt\\
\nonumber
	&\quad - \int_{\mathbb{R}^n_+} ( (u_0(x)+u_1(x))\psi(0,x) - u_0(x) \psi_t (0,x))\,dx,
\end{align}%
where we used the notation
$x' = (x_1, \ldots, x_{n-1})$.
Now, we choose the test function $\psi$ as follows.
Let $\eta (t) \in C_0^{\infty}([0,\infty))$ be a non-increasing function satisfying
\begin{align*}%
	\eta(t) = 1 \ (t \in [0,1/2]), \quad \eta (t)= 0 \ (t \in [1, \infty)).
\end{align*}%
We also define
$\phi \in C_0^{\infty}(\mathbb{R}^n)$ by
$\phi(x) := \eta(|x_1|) \eta(|x_2|) \cdots \eta(|x_n|)$.
Let $R> 0$ be a parameter and let
$\psi_R(t,x) := \phi( x/ R) \eta (t/R^2)$.
We denote the rectangle
$D_R := \{ x \in \mathbb{R}^n \,;\, |x_1| \le R, \ldots, |x_n|\le R \}$
and we put $D_R^+ = D_R \cap \mathbb{R}^n_+$.
Then, it is obvious that
${\rm supp\,} ( \partial_{x_j} \phi (\cdot /R) ) \subset D_R \setminus D_{R/2}$.
With the above notations, we choose our test function as
$\psi(t,x) = x_n \psi_R(t,x)^l$ with sufficiently large integer $l$.

Let
\begin{align*}%
	I_R := \int_0^{\infty}\int_{\mathbb{R}^n_+} |u|^p x_n \psi_R^l \,dxdt.
\end{align*}%
We note that our choice of test function implies
$\displaystyle \left. x_n \psi_R(t,x)^l \right|_{x_n=0} = 0$
and
$\displaystyle \left. x_n \partial_t (\psi_R(t,x)^l) \right|_{t=0} = 0$.
Moreover, by the assumption \eqref{ini}, we see that
there exists $R_0>0$ such that
\begin{align*}%
	\int_{\mathbb{R}^n_+} ( (u_0(x)+u_1(x)) x_n \psi_R(0,x)^l\, dx > 0
\end{align*}%
holds for $R \ge R_0$.
Therefore, we deduce form \eqref{int} that
\begin{align*}%
	I_R &\le \int_0^{\infty}\int_{\mathbb{R}^n_+}
			u (\partial_t^2 (x_n \psi_R^l) - \Delta (x_n \psi_R^l) - \partial_t(x_n \psi_R^l) ) \,dxdt\\
		&=: K_1 + K_2 + K_3
\end{align*}%
for $R \ge R_0$.
We estimate $K_1, K_2$ and $K_3$ individually.
First, for $K_1$,
we apply the H\"{o}lder inequality to obtain
\begin{align*}%
	K_1 &\le CR^{-4}
		\left( \int_{R^2/2}^{R^2}\int_{D_R^+} |u|^p x_n \psi_R^l \,dxdt\right)^{1/p}
		\left( \int_{R^2/2}^{R^2}\int_{D_R^+} x_n \,dxdt\right)^{1/p'}\\
	&\le CR^{-4+(n+3)/p'} \hat{I}_R^{1/p},
\end{align*}%
where
$p'$ stands for the H\"{o}lder conjugate of $p$ and
\begin{align*}%
	\hat{I}_R := \int_{R^2/2}^{R^2}\int_{D_R^+} |u|^p x_n \psi_R^l \,dxdt.
\end{align*}%
Similarly, by using
\begin{align*}%
	\Delta (x_n \psi_R^l)
	&= l R^{-2} x_n \left( \phi\left( \frac{x}{R} \right)^{l-1} (\Delta \phi) \left( \frac{x}{R} \right)
		+ (l-1) \phi\left( \frac{x}{R} \right)^{l-2}
				\left| (\nabla \phi) \left( \frac{x}{R} \right) \right|^2 \right)
				\eta \left( \frac{t}{R^2} \right)^l\\
	&\quad + 2l R^{-1} \phi \left( \frac{x}{R} \right)^{l-1}
			(\partial_{x_n} \phi) \left( \frac{x}{R} \right) 
				\eta \left( \frac{t}{R^2} \right)^l,
\end{align*}%
we estimate $K_2$ as
\begin{align*}%
	K_2 &\le 
	C R^{-2} \left(
		\int_0^{R^2} \int_{D_R^+\setminus D_{R/2}^+}|u|^p x_n \psi_R^l \,dxdt
			\right)^{1/p}
		\left( \int_0^{R^2} \int_{D_R^+\setminus D_{R/2}^+} x_n \,dxdt \right)^{1/p'}\\
	&\quad + CR^{-1} 
		\left(
		\int_0^{R^2} \int_{D_R^+\setminus D_{R/2}^+}|u|^p x_n \psi_R^l \,dxdt
			\right)^{1/p}
		\left( \int_0^{R^2} \int_{D_R^+ \cap \{ x_n > R/2\} } x_n^{-p'/p} \,dxdt \right)^{1/p'}\\
	&\le CR^{-2 + (n+3)/p'} \tilde{I}_R^{1/p},
\end{align*}%
where
\begin{align*}%
	\tilde{I}_R
	= \int_0^{R^2} \int_{D_R^+\setminus D_{R/2}^+}|u|^p x_n \psi_R^l \,dxdt
\end{align*}%
and we note that $(\partial_{x_n}\phi)(x/R) = 0$ on the set $\{ x_n \le R/2\}$.
The term $K_3$ can be estimated in the same way as $K_1$ and we have
\begin{align*}%
	K_3 \le CR^{-2+(n+3)/p'} \hat{I}_R^{1/p}.
\end{align*}%
Combining the estimates above, we deduce
\begin{align}%
\label{IR}
	I_R \le C( R^{-4+(n+3)/p'} \hat{I}_R^{1/p}
		+ R^{-2 + (n+3)/p'} \tilde{I}_R^{1/p}
		+ R^{-2+(n+3)/p'} \hat{I}_R^{1/p}).
\end{align}%
In particular, using
$\hat{I}_R \le I_R$ and $\tilde{I}_R \le I_R$,
we have
\begin{align}%
\label{IR2}
	I_R \le C ( R^{-4+(n+3)/p'} + R^{-2 + (n+3)/p'}) I_R^{1/p}.
\end{align}%
When $1<p<p_F(n+1)$, letting $R \to \infty$, we see that
$I_R \to 0$, which implies $u \equiv 0$.
However, this contradicts $(u_0, u_1) \not\equiv 0$.

On the other hand, when $p = p_F(n+1)$,
we have $-2+(n+3)/p'=0$ and hence,
we see from \eqref{IR2} that
$I_R \le C$ with a constant $C$ independent of $R$.
Thus, letting $R\to \infty$, we have
$x_n |u|^p \in L^1([0,\infty) \times \mathbb{R}^n_+)$.
Noting this and the integral region of
$\hat{I}_R$ and $\tilde{I}_R$,
we also deduce
\begin{align*}%
	\lim_{R\to \infty} ( \hat{I}_R + \tilde{I}_R ) = 0.
\end{align*}%
This and \eqref{IR} imply
\begin{align*}%
	I_R \le C ( \hat{I}_R^{1/p} + \tilde{I}_R^{1/p} ) \to 0\quad (R\to \infty),
\end{align*}%
and hence, $u \equiv 0$.
This again contradicts $(u_0, u_1) \not\equiv 0$
and completes the proof.


\begin{thebibliography}{9}
\bibitem{Fu66}\textsc{H. Fujita},
{\em On the blowing up of solutions of the Cauchy problem for $u_t=\Delta u+u^{1+\alpha}$},
J. Fac. Sci. Univ. Tokyo Sec. I, {\bf 13} (1966), 109-124.

\bibitem{GeYaLa16}{\sc J. Geng, Z. Yang, N. Lai},
{\em Blow-up and lifespan estimates for initial boundary value problems
for semilinear Schr\"{o}dinger equations on half-line} (Chinese),
Acta Math. Sci. Ser. A Chin. Ed. {\bf 36} (2016), 1186--1195. 

\bibitem{Hay73}{\sc K. Hayakawa},
{\em On nonexistence of global solutions of some semilinear parabolic differential equations},
Proc.\ Japan Acad.\ {\bf 49} (1973), 503--505.

\bibitem{HaKaNa04JAA}{\sc N. Hayashi, E. I. Kaikina, P. I. Naumkin},
{\em Damped wave equation with a critical nonlinearity on a half line},
J.\ Anal.\ Appl.\ {\bf 2} (2004), 95--112.

\bibitem{Ik03DIE}{\sc R. Ikehata},
{\em A remark on a critical exponent for 
the semilinear dissipative wave equation in the one dimensional half space},
Differential Integral Equations {\bf 16} (2003), 727--736.

\bibitem{Ik03JMAA}{\sc R. Ikehata},
{\em Critical exponent for semilinear damped wave equations in the $N$-dimensional half space},
J. Math. Anal. Appl. {\bf 288} (2003), 803--818. 

\bibitem{Ik04MMAS}{\sc R. Ikehata},
{\em New decay estimates for linear damped wave equations
and its application to nonlinear problem},
Math. Methods Appl. Sci. {\bf 27} (2004), 865--889.

\bibitem{KoShTa77}{\sc K. Kobayashi, T. Shirao, H. Tanaka},
{\em On the growing up problem for semilinear heat equations},
J. Math. Soc. Japan {\bf 29} (1977), 407--424. 

\bibitem{LeMe89}{\sc H.A. Levine, P. Meier},
{\em A blowup result for the critical exponent in cones},
Israel J. Math. {\bf 67} (1989) 129--136.

\bibitem{LeMe90}{\sc H.A. Levine, P. Meier},
{\em The values of the critical exponent for reaction-diffusion equation in cones},
Arch. Rational Mech. Anal. {\bf 109} (1990) 73--80.


\bibitem{Me88}{\sc P. Meier},
{\em Blow up of solutions of semilinear parabolic differential equations},
J. Appl. Math. Phys. {\bf 39} (1988) 135--149.

\bibitem{Me90}{\sc P. Meier},
{\em On the critical exponent for reaction-diffusion equations},
Arch. Rational Mech. Anal. {\bf 109} (1990) 63--71.

\bibitem{NiZh05}{\sc K. Nishihara, H. Zhao},
{\em Existence and nonexistence of time-global solutions to damped wave equation on half-line},
Nonlinear Anal. {\bf 61} (2005), 931--960.

\bibitem{ToYo01}{\sc G. Todorova, B. Yordanov},
{\em Critical exponent for a nonlinear wave equation with damping},
J. Differential Equations {\bf 174} (2001), 464--489.

\bibitem{Zh01}{\sc Qi S. Zhang},
{\em A blow-up result for a nonlinear wave equation with damping: the critical case},
C. R. Acad. Sci. Paris S\'{e}r. I Math. {\bf 333} (2001), 109--114.

\end{thebibliography}
\end{document}